\documentclass{siamltex}
\usepackage{amstext,amsmath,amssymb}         
\usepackage{mathtext}
\usepackage{multicol}
\usepackage{pgfplots}
\usepackage{makecell}
\usepackage{placeins}
\usepackage{multicol}
\usepackage{pgfplots}
\usepackage[colorinlistoftodos,textwidth=4cm,shadow]{todonotes}
\usepackage[pdftex]{hyperref}
\usepackage{algorithmic,algorithm}

\newcounter{Ivan}

\newcounter{Maxim}

\pdfoutput=1
\title{\bf Fast multidimensional convolution in low-rank tensor formats via cross approximation}

\author{M.~V. Rakhuba\footnotemark[2] \footnotemark[3]
\and I.~V. Oseledets\footnotemark[2] \footnotemark[4]}

\begin{document}
\maketitle
\renewcommand{\thefootnote}{\fnsymbol{footnote}}

\footnotetext[2]{Skolkovo Institute of Science and Technology, Novaya St. 100, 143025 Skolkovo, Moscow Region, Russia.}
\footnotetext[3]{Moscow Institute of Physics and Technology, Institutskii per. 9, Dolgoprudny, 141700, Moscow Region, Russia}
\footnotetext[4]{Institute of Numerical Mathematics, Russian Academy of Sciences. Gubkina St. 8, 119333 Moscow, Russia.}
\renewcommand{\thefootnote}{\arabic{footnote}}




\begin{abstract}
We propose new {\it cross-conv} algorithm for approximate computation of convolution in different low-rank tensor formats (tensor train, Tucker, Hierarchical Tucker).
It has better complexity with respect to the tensor rank than previous approaches.
The new algorithm has a high potential impact in different applications.
The key idea is based on applying cross approximation in the ``frequency domain'', where convolution becomes a simple elementwise product.
We illustrate efficiency of our algorithm by computing the three-dimensional Newton potential 
and by presenting preliminary results for solution of the Hartree-Fock equation on tensor-product grids.
\end{abstract}



\section{Introduction}
Multivariate convolution problem arises in a range of applications, such as population balance models \cite{cor-pbm-2014}, Smoluchowski equation \cite{wattis-smoluchowskieq-2006,palaniswaamy-smol-2007}, modeling of quantum mechanical systems with the help of Hartree-Fock and Kohn-Sham equations \cite{beylkin-quantum-2004,beylkin-multires-exchange-2004,vekh-hartree-2008,venera-phd}. 
Several applications can be found in signal/data processing \cite{woods-multisig-2006} or even in financial mathematics \cite{kwok-option-2012}. 

The convolution of $f,g :\mathbb{R}^d\rightarrow \mathbb{R} $ is defined by the integral transform 
\begin{equation}\label{cont-conv}
	(f*g)(x) \equiv \int_{\mathbb{R}^d}  f(y)\, g(x-y) \, dy, \quad x\in \mathbb{R}^d,
\end{equation}
where $f$ is assumed to have bounded support.
We get a discrete convolution problem by a suitable discretization of \eqref{cont-conv} on a uniform grid:
\begin{equation}\label{discr-conv}
	(f*g)_\bold{i} = \sum_{\bold{j}} f_\bold{j}\, g_{\bold{i}-\bold{j}}, 
\end{equation} 
where $\bold{i},\bold{j}\in \{0,\dots, n-1 \}^d$ are multi-indices. Usage of uniform grids is typical, but non-obligatory. 
Non-uniform grids \cite{hackbush-nonequid-2007, hackbush-nonequid-2008} can be used. In this paper we consider only uniform grids and  the  discrete convolution \eqref{discr-conv} is the main object of study in this paper.

Classic approach to compute the discrete convolution is based on the Fast Fourier Transform (FFT). 
It requires $\mathcal{O}(n^d \log n)$  operations for a grid with $n^d$ points.  This is much faster than the naive approach (with complexity $\mathcal{O}(n^{2d})$), but 
still prohibitive for large $d$ and/or $n$. To reduce computational complexity certain low-parametric representations of $f$ and $g$ have to be used.
For this task we will use \emph{tensor formats} which are based on the idea of \emph{separation of variables}. 
The most straightforward way to separate variables is to use the canonical polyadic format (CP format, also called CANDECOMP/PARAFAC model) which dates back to 1927 \cite{hitchcock-rank-1927}. 
A tensor is said to be in canonical format if it can be represented in the form
$$
	A(i_1, \ldots, i_d) = \sum_{\alpha=1}^{r} U_1(i_1, \alpha) \dots U_d(i_d, \alpha),
$$
where the minimal possible $r$  is called \emph{canonical rank}. If a good CP approximation is known,
many basic operations are fast to compute \cite{beylkin-2002, khor-rstruct-2006, khor-low-rank-kron-P1-2006,khor-low-rank-kron-P2-2006,khor-tuckertype-2007,ost-latensor-2009,beylkin-high-2005}. 

Nevertheless, the CP decomposition suffers from a serious drawback: there are no robust algorithms to compute it numerically  for $d>2$ \cite{desilva-2008}. 
Note that in two dimensions it can be computed in a stable way by using SVD or, if the matrix is large, by rank-revealing algorithms.

The Tucker format \cite{Tucker, lathauwer-svd-2000,lathauwer-rank1-2000,khor-tuckertype-2007} is another classic decomposition of tensors. 
It can be computed via stable algorithms but the number of parameters grows exponentially in $d$. 
As a result, it is typically used only for problems with small $d$, especially for the three-dimensional case.
In higher dimensions other \emph{stable tensor formats}, namely tensor train (TT) \cite{ot-tt-2009, osel-tt-2011} or  hierarchical Tucker (HT) \cite{hk-ht-2009, gras-hsvd-2010} formats can be used. In contrast with the Tucker format, they do not suffer from the ``curse of dimensionality''. For more details regarding low-rank representations of tensors see the book by Hackbusch \cite{hackbusch-2012} and reviews \cite{khor-survey-2011, larskres-survey-2013, kolda-review-2009}.  

{ \it Related work.} In this paper we focus on fast computation of multidimensional convolution. 
Although it is not difficult to implement convolution in complexity linear in $d$ or $n$, a \emph{strong rank dependence may occur}. 
The rank of the result is generally equal to the product of the ranks of $f$, $g$, and then one should \emph{truncate} the representation with necessary accuracy (by truncation we mean approximation in the same format with smaller rank).
This approach was considered in \cite{st-chem-2009, khor-acc-2010} and may lead to high complexity when the ranks are large.
A remarkable work is \cite{khkaz-conv-2013} where an algorithm for the computation of convolution in so-called Quantized TT (QTT) \cite{khor-qtt-2011, osel-2d2d-2010} was proposed. This algorithm has complexity $\mathcal{O}(d \log^{\alpha} n)$ and is asymptotically the best one. However, for $n$ of practical interest the algorithm proposed in this paper is faster for the same discretization and approximation accuracy $\varepsilon$. This is due to high constant hidden in $\mathcal{O}(\cdot)$ term in the QTT algorithm.

The algorithm proposed in this paper is simple.
At first, we use a classic idea of representing discrete convolution in the form of several Fourier transforms and one element-wise multiplication in the ``frequency domain''. 
The crucial step is to interpolate this element-wise product via \emph{cross approximation method}. One of the nice properties of the stable (SVD-based) tensor formats (Tucker, TT, HT) is that for each of them there is an algorithm, that allows to accurately reconstruct a low-rank tensor using only few of its elements. The Fourier transform steps do not change the tensor structure,  and the element-wise multiplication is done via the cross approximation algorithm. 

Our paper is organized as follows. 
In Section \ref{sec-notations} we give a brief summary of notations.
In Section \ref{sec-discretization} we discuss different discretizations that lead to the discrete convolution.
The cross-conv algorithm is described in Section \ref{sec-algorithm} and its complexity is analyzed.
In Section \ref{sec-experiments} numerical experiments are  presented: we compute three-dimensional Newton potentials of different electronic densities. 
We also compare our algorithm with one in \cite{khkaz-conv-2013}. 
Finally, we present preliminary results for the solution of the Hartree-Fock equation on tensor-product grids.
In Appendix \ref{sec-cross3d} we present new cross approximation algorithm -- Schur-Cross3D.

\section{Notation and prerequisites}\label{sec-notations}
In this section we will give a brief summary of notations that we use. This material is not new and 
can be found in \cite{kolda-review-2009,larskres-survey-2013,hackbusch-2012}.

{\it Tensors} are just multidimensional arrays. They will be denoted by boldface letters, i.e. $\bold{A}$. We denote an element of $\bold{A}$ in position $(i_1,\ldots,i_d)$ as  $A(i_1, \dots, i_d)$. The number of indices $d$ will be called \emph{dimension} of a tensor. 
 Indices $i_k$ vary from $0$ to $n_k-1$  (this makes the notation for the convolution simpler), where $n_k$ are called \emph{mode sizes}.
The Frobenius norm of a tensor is defined as
\[ 
	\Vert \bold{A} \Vert = \sqrt{\sum_{i_1,\dots, i_d} |A(i_1,  \dots, i_d)|^2}. 
\]
The element-wise tensor multiplication of tensors $\bold{A}$ and $\bold{B}$ is denoted by $\bold{C}=\bold{A}\circ \bold{B}$ and is defined as
\[
	C({i_1, \dots, i_d}) = A(i_1, \dots, i_d) B(i_1, \dots, i_d).
\]
A tensor $\bold{A}$ is said to be in the Tucker format \cite{Tucker}, if it is represented as
\begin{equation} \label{Tucker}
	A(i_1, \dots, i_d) = \sum_{\alpha_1, \dots, \alpha_d} G^{(\bold{A})} (\alpha_1, \dots, \alpha_d) \, U_1^{(\bold{A})}(i_1, \alpha_1) \dots U_d^{(\bold{A})}(i_d, \alpha_d),
\end{equation}
where $\alpha_k$ varies from $1$ to $r_k$. 
The minimal number of summands $r_k$ required to represent $\bold{A}$ in the form \eqref{Tucker} is called
the Tucker rank of $k$-th mode. The tensor $\bold{G}^{(\bold{A})} $ is called the core of the decomposition and $U_k^{(\bold{A})} $ are referred to as Tucker factors.
The Tucker decomposition contains $\mathcal{O}(r^d + nrd)$ elements, so the number of parameters grows exponentially in $d$. 

Tensor train (TT) (or MPS in other communities) and Hierarchical Tucker (HT) formats are efficient low-parametric representations of multidimensional tensors. 
A tensor $\bold{A}$ is said to be in the TT-format \cite{ot-tt-2009,osel-tt-2011}  if it can be written in the form
\begin{equation} \label{tt}
	A(i_1, \dots, i_d) = \sum_{\alpha_0, \dots, \alpha_d} G_{1}^{(\bold{A})} (\alpha_0,i_1, \alpha_1) \, G_{2}^{(\bold{A})} (\alpha_1, i_2, \alpha_2) \dots  G_{d}^{(\bold{A})} (\alpha_{d-1}, i_d, \alpha_d).
\end{equation}
In \eqref{tt} $G_k$ have sizes $r_{k-1} \times n_k \times r_k$ and are called \emph{TT-cores}, where $r_0=1$ and $r_d = 1$.
The numbers $r_k$ are called TT-ranks of the representation. The decomposition \eqref{tt} can be also written in the matrix-product form (Matrix Product State, MPS)
\[
	A(i_1, \dots, i_d) = G_1^{(\bold{A})} (i_1) G_2^{(\bold{A})}(i_2) \dots G_d^{(\bold{A})}(i_d),
\]
where $G_k(i_k)$ are $r_{k-1}\times r_k$ matrices that depend on parameter $i_k$.
It is worth to note that the MPS representation, which is algebraically equivalent to the TT-format, has been used for a long time in quantum information theory and solid state physics to approximate certain wavefunctions \cite{white-dmrg-1992, ostlund-dmrg-1995}, see the review \cite{schollwock-2011} for more details.

An alternative way to reduce the complexity in the multidimensional case is given by the
hierarchical Tucker (HT) decomposition.
The idea is to apply Tucker decomposition recursively by merging indices according to some binary-dimension tree. 
The tensor is represented by a collection of \emph{transfer tensors} corresponding to the nodes of the tree. For the  linear tree HT-format reduces to the TT-format. In practice, however, much simpler structure of the TT-format is more convenient for the implementation of different algorithms.

Discrete Fourier Transform is crucially required for fast convolution algorithms. We denote by $\mathcal{F} (\bold{A})$ Fourier transform of the tensor $\bold{A}$:
$$\mathcal{F} (\bold{A})(i_1, \dots, i_d) = \sum_{j_1,\dots, j_d} e^{-2\pi i \left[ \frac{i_1j_1}{n_1} + \dots + \frac{i_dj_d}{n_d}\right]} A(j_1, \dots, j_d),$$
and by $\mathcal{F}^{-1} (\bold{A})$ inverse Fourier transform:
$$\mathcal{F}^{-1} (\bold{A})(i_1, \dots, i_d) = \frac{1}{n_1\dots n_d}\sum_{j_1,\dots, j_d} e^{2\pi i \left[ \frac{i_1j_1}{n_1} + \dots + \frac{i_d j_d}{n_d}\right]} A(j_1, \dots, j_d).$$

\section{Discretization}\label{sec-discretization}
For convenience we describe here well-known facts about the discretization of the convolution, see, for example, \cite{khor-acc-2010}.
Recall that the multidimensional convolution of $f,g :\mathbb{R}^d\rightarrow \mathbb{R} $ is defined by the integral transform 
\[ (f*g)(x) \equiv \int_{\mathbb{R}^d} f(y) g(x-y) \, dy, \quad x\in \mathbb{R}^d. \]
We assume that the convolving function $f$ is from $\mathcal{L}_2(\mathbb{R}^d)$ and has bounded support in the box $\Omega = [-L, L]^d$.
The size of the box depends on the application, but in many cases (i.e. in electronic structure computation) functions decay exponentially with $\Vert x \Vert \rightarrow \infty$ and the choice is obvious.
The function $g$ is such that the $f*g$ is from $\mathcal{L}_2(\mathbb{R}^d)$. 
In particular we are interested in the calculation of the Newton potential where $g(x) = 1/\|x\|$.
Note that in the general case the convolution is not a continuous mapping from $\mathcal{L}_2\times \mathcal{L}_2$ into $\mathcal{L}_2$.
However, once the convolution is discretized, the convergence in one norm means convergence in any other due to norm equivalence.
The estimate for the convolution in the spectral norm is presented in Section \ref{sec-accuracy}.

There are three standard ways to discretize convolution: Galerkin method, collocation method and Nystr\"om-type schemes. 
First, introduce in $\Omega$ a uniform  tensor-product grid $\omega^h = \omega_1^h \times \dots \times \omega_d^h$ with $h = 2L/n$, where $\omega_i^h = \{ -L + kh: \, k=0, ..., n \}$, $i=1,\dots,d$. 
For simplicity, consider piecewise-constant basis functions $\phi_{\bold{i}}$ with support on $\Omega_{\bold{i}}$, where $\bold{i} \in \mathcal{I} \equiv \{0,\dots, n-1\}^d$ and $\Omega_{\bold{i}}$ are cubes with edge size $h$ centered in $y_{\bold{i}}$. 
Thus, we have
\begin{equation}\label{1}
  (f*g)(x) \approx \sum_{\bold{i}\in \mathcal{I}} f_{\bold{i}}  \int_{\Omega_{\bold{i}} } \phi_{\bold{i}}(y)  g(x-y) \, dy,
\end{equation}
where $f_{\bold{i}}$ are the coefficients in the expansion  $f(y) \approx \sum_{\bold{i}\in\mathcal{I}} f_{\bold{i}} \phi_{\bold{i}}(y)$.
As a result, the collocation scheme with collocation points $x_{\bold{j}}$ located on uniform tensor-product grid with the same same step size $h$ yields a discrete convolution:
\begin{equation}\label{coll}
w_{\bold{j}} \equiv (f*g)(x_{\bold{j}})  \approx \sum_{\bold{i}} f_\bold{i}\, g_{\bold{i}-\bold{j}}   ,\quad \bold{j} \in \mathcal{I},
\end{equation}
where 
\begin{equation} \label{coef}
	g_{\bold{i}-\bold{j}} =  \int_{\Omega_{\bold{i}} } \phi_{\bold{i}}(y)  g(x_{\bold{j}}-y) \, dy,
\end{equation}
is a multilevel Toeplitz matrix. The problem with the collocation method is that it leads to non-symmetric Toeplitz matrices even if the original convolution 
was symmetric. This may pose problems in some applications. A natural choice is to use a Galerkin method,
which again leads to the discrete convolution with
\begin{equation} \label{galerk}
	g_{\bold{i}-\bold{j}} =  \int_{\mathbb{R}^d} \phi_{\bold{i}}(x) \phi_{\bold{j}}(y)  g(x-y) \, dx dy, \quad f_{\bold{i}} = \int_{\mathbb{R}^d} f(x) \phi_{\bold{i}}(x)  \, dx.
\end{equation}
To get high-order discretization schemes translation-invariant basis-functions of higher order can be used $\phi_\bold{i}(y) = \psi(y - y_\bold{i})$, 
where $\psi(y)$ is a suitable piecewise-polynomial function. Computation of matrix elements in \eqref{coef} or \eqref{galerk} even for piecewise-constant functions can be complicated. A simple alternative is a Nystr\"om-type scheme that uses shifted grids \cite{dks-ttfft-2012}
\begin{equation}\label{nystrom}
  (f*g)(x_\bold{j}) \approx h^d \sum_{\bold{i}\in \mathcal{I}} f(y_{\bold{i}})    g(x_\bold{j}-y_\bold{i}), \quad \bold{j} \in \mathcal{I},
\end{equation}
where $x_\bold{j}$ are points of $y_\bold{i}$ shifted by half step. For a certain class of functions it provides almost second order of accuracy up to a logarithmic term.


\section{Algorithm description}\label{sec-algorithm}
Let us consider  a $d$-dimensional discrete convolution of two tensors $f_{\bold{i}}$ and $g_{\bold{j}}$
\begin{equation}\label{conv}
w_{\bold{j}} = \sum_{\bold{i}\in \mathcal{I}}    f_\bold{i}\, g_{\bold{i}-\bold{j}} ,\quad \bold{j} \in \mathcal{I}.
\end{equation}
This can be also considered as a product of a \emph{multilevel Toeplitz matrix}  with elements $g_{\bold{i}-\bold{j}}$ by a vector (see, for example, \cite{tee-circopt-1992} for 
properties of multilevel Toeplitz matrices). The computation of \eqref{conv} as a direct sum requires $\mathcal{O}(n^{2d})$ operations.  Using the FFT the complexity can be reduced to $\mathcal{O}(n^{d} \log n)$.  The classic FFT-based algorithm is our starting point for an efficient low-rank convolution algorithm.

The idea of  the FFT method is to replace the Toeplitz matrix by vector product to the product of a larger \emph{circulant} matrix  by vector. 
For instance, a 1-level $n\times n$ Toeplitz matrix $ \{g_{i-j} \}_{i,j =0}^{n-1}$ may be embedded in an $(2n-1)\times(2n-1)$ circulant matrix which is fully defined by its 
first column $\bold{c}_g \equiv \{g_0, g_1, \dots, g_{n-1}, g_{1-n}, g_{2-n}, \dots, g_{-1}\}$.

In the $d$-dimensional case a multilevel circulant matrix is defined by a tensor $\bold{c}_g$:
$$
c_g(i_1, \dots, i_d)= g_{\tau({i_1}),  \dots, \tau({i_d})}, \quad i_1,\dots,i_d\in \overline{0, \, 2n-2},
$$ 
where
$$
\tau(i)=
\begin{cases} 
	i,   &i\in \overline{0, \, n-1},\\  
	i-2n+1,  &i \in \overline{n,\, 2n  -2}.    
\end{cases}
$$ 
At the first step we embed $\bold{f}$ into a larger tensor $\bold{q}_f$ with mode sizes $(2n_1-1, \ldots, 2 n_d -1)$ by zero-padding:
$$
q_f(i_1, \dots, i_d)=
\begin{cases} 
	f_{i_1, \dots, i_d}, & \quad i_1,\dots,i_d\in \overline{0,\, n-1},\\  
	0, & \quad \text{otherwise.}
\end{cases}
$$ 
Multilevel circulant matrices are diagonalized by the normalized unitary Fourier matrix $1/n^{\frac{d}{2}}\ F_d$ and the eigenvalues can be computed from the DFT of the first column, 
$$
    C = \frac{1}{n^d} F_d^* \Lambda F_d,
$$
where
$$\Lambda  = \textrm{diag}(\mathcal{F}(\bold{c}_g)).$$
Therefore,
\begin{equation}\label{3}
\bold{\tilde w} =  \mathcal{F}^{-1}\left( \mathcal{F} (\bold{c}_g)  \circ \mathcal{F} (\bold{q}_f) \right),
\end{equation}
where $\bold{\tilde w}$ 
is the expanded convolution tensor with $(2n-1)$ each mode size and we are interested only in its subtensor $\bold{w}$:
$$
	w(i_1, \dots, i_d) = \tilde w (i_1, \dots, i_d), \quad i_1,\dots,i_d\in \overline{0, \, n-1}.
$$
How to use this formula if the operands are given in a low-rank tensor format? For simplicity, consider that 
$\bold{c}_g$ and $\bold{q}_f$ are in the TT-format
\begin{equation}\label{cq-tt}
\begin{aligned}
 c_g( i_1, \dots, i_d) &= G^{(\bold{c}_g)}_1(i_1) \dots  G^{(\bold{c}_g)}_d(i_d), \\
  q_f ( i_1, \dots, i_d) &= G^{(\bold{q}_f)}_1(i_1) \dots  G^{(\bold{q}_f)}_d(i_d),
\end{aligned}
\end{equation}
however the idea applies to other SVD-based formats (Tucker, HT, skeleton).
The Fourier matrix has a tensor product structure:
$$
  F_d = F\otimes F \otimes \ldots \otimes F,
$$
therefore, its application does not change the TT-ranks (as well as the inverse Fourier transform).
Indeed, given a tensor $\bold{A}$ in the TT-format, the $\mathcal{F}(\bold{A})$ can be written in the following form:
\begin{equation} \label{fft-sep}
\begin{aligned}
\mathcal{F} (\bold{A})(i_1, \dots, i_d) &= \sum_{j_1,\dots, j_d} e^{-2\pi i \left[ \frac{i_1j_1}{n_1} + \dots + \frac{i_d j_d}{n_d}\right]} G_1^{(\bold{A})} (j_1) \dots G_d^{(\bold{A})}(j_d) = \\
&= \sum_{j_1} e^{-2\pi i  \frac{i_1j_1}{n_1}} \, G_1^{(\bold{A})} (j_1)  \dots \sum_{j_d} e^{-2\pi i  \frac{i_dj_d}{n_d}}\, G_d^{(\bold{A})} (j_d) = \\
&=\mathcal{F}_{1D} \left(G_1^{(\bold{A})}\right) (i_1) \dots \mathcal{F}_{1D} \left(G_d^{(\bold{A})}\right) (i_d), 
\end{aligned}
\end{equation}
where by $\mathcal{F}_{1D}$ we denote a 1-dimensional Fourier transform.

Now we are ready to describe the algorithm.

\subparagraph{Step 1} Compute tensors $\mathcal{F} (\bold{c}_g)$ and $\mathcal{F} (\bold{q}_f)$ in the considered format. 
As was mentioned above, Fourier transform of any tensor does not change its ranks and is equivalent to univariate Fourier transforms of each factor in the Tucker case 
and each core in the TT case. 
Therefore in the TT-format,
\begin{equation}\label{cq-tt1}
\begin{aligned}
 \mathcal{F}(\bold{c}_g)( i_1, \dots, i_d) &= \mathcal{F}_{1D} \left(G_1^{(\bold{c}_g)}\right) (i_1)  \dots\mathcal{F}_{1D} \left(G_d^{(\bold{c}_g)}\right) (i_d), \\
 \mathcal{F}(\bold{q}_f)( i_1, \dots, i_d) &= \mathcal{F}_{1D} \left(G_1^{(\bold{q}_f)}\right) (i_1)  \dots\mathcal{F}_{1D} \left(G_d^{(\bold{q}_f)}\right) (i_d),
\end{aligned}
\end{equation}

\subparagraph{Step 2}
In this step we compute element-wise product $\bold{\Theta} = \mathcal{F} (\bold{c}_g)  \circ \mathcal{F} (\bold{q}_f)$ and this is the crucial step of our algorithm. 
The naive approach is to compute it directly and it leads to the tensor representation with the ranks squared. The truncation is almost always required, and for the  TT-format it
leads to an algorithm with complexity $\mathcal{O}(dnR^3)$, where $R = r^2$. Such algorithm works only up to ranks $r_k \sim 100$. There are more sophisticated algorithms for different formats that are based on iterative schemes, e.g. for Tucker \cite{gos-kryl-2012, stz-tr-2012} and for TT-format \cite{Os-mvk2-2011, ds-amen-2014}, but they work on the resulting tensor.
We propose to compute the element-wise product via \emph{sampling}.  It is very cheap to compute any prescribed element of the product,
and that situation is perfectly suited for the application of \emph{cross approximation methods}. Such methods are proposed for all of the SVD-based formats and require the same amount
of elements to be sampled, as the number of parameters in the decomposition! We will give corresponding complexity estimates in the next section. Thus, we compute necessary elements of the tensors $\mathcal{F} (\bold{c}_g)$ and $\mathcal{F} (\bold{q}_f)$, multiply them and build a tensor $\bold{\Theta}$ in the considered format according to those elements and selected cross approximation scheme. 
This is the only step where approximation is done. Suppose that the approximation error is $\delta$:
\[
{\Theta}( i_1, \dots, i_d)  = \widetilde{\Theta} + \Delta \bold{\Theta},	
\]
where
\[
{\widetilde\Theta}( i_1, \dots, i_d)  = G_1^{(\bold{\Theta})} (i_1) \dots G_d^{(\bold{\Theta})} (i_d),	
\]
is the approximation of $\Theta$ computed via a cross method with relative accuracy \\
$ \|\Delta \bold{\Theta} \|/\|\bold{\Theta}\| = \delta$.

\subparagraph{Step 3}
Compute $\mathcal{F}_{1D}^{-1}$ of each $\bold{\widetilde{\Theta}} $ core. 
Therefore, the final approximation $\bold{\tilde w}$ has the form
\[
	\bold{\tilde w} = \mathcal{F}^{-1} (\bold{\Theta})( i_1, \dots, i_d)  = \mathcal{F}^{-1}_{1D} \left(G_1^{(\bold{\Theta})}\right) (i_1) \dots \mathcal{F}^{-1}_{1D} \left(G_d^{(\bold{\Theta})}\right) (i_d) + \mathcal{F}^{-1} (\Delta \Theta) .
\]
It is easy to estimate the required threshold $\delta$ to be provided to the cross approximation algorithm. Suppose $\epsilon = \|\Delta \bold{\tilde w} \|/\|\bold{\tilde w}\|$, is the required accuracy, where $\Delta \bold{\tilde w} = \mathcal{F}^{-1} (\Delta \Theta)$. 
Due to unitary invariance of the Frobenius norm we have
\[
	\epsilon = \frac{\|\Delta \bold{\tilde w} \|}{\|\bold{\tilde w}\|} = 
			  \frac{\|\mathcal{F}^{-1} (\Delta \Theta)\|}{\|\mathcal{F}^{-1} (\Theta)\|} = 
	                 \frac{\|\Delta \Theta\|}{\|\Theta\|} = \delta.
\]	
So, to provide the convolution accuracy $\epsilon$ one needs to run the cross approximation algorithm with $\delta(\epsilon) = \epsilon$.

\subsection{Controlling the accuracy for perturbed inputs} \label{sec-accuracy}
Note that there are two sources of errors. 
First source is a discretization scheme error. 
The only interesting point for us is that the chosen discretization approximates the convolution with the required order. 
Once the continuous convolution is reduced to the discrete one, another source of the error is connected with low-rank approximation. 
We control this error by using the euclidean of vectors which may be not the optimal choice on the discrete level.
Indeed, let us consider convolution as a matrix-by-vector multiplication $w = G f$, where $G$ is a Toeplitz matrix generated by the vector $g$. 
Let us estimate $\Delta w$ in $w+\Delta w = (G + \Delta G) (f + \Delta f)$ where $\Delta G$ and $\Delta f$ are small pertrubations connected with tensor approximations. Up to second order corrections
$$
\| \Delta w \|_2 \leqslant \|\Delta G\|_2 \|f\|_2 + \|G\|_2 \|\Delta f\|_2.
$$
Since $G$ and $\Delta G$ are Toeplitz and considering the fact that $\|\cdot\|_2 \leqslant \|\cdot\|_F$ we have
$$
\| \Delta w \|_2 \leqslant n^{d/2} \left(\|\Delta g\|_2 \|f\|_2 + \|g\|_2 \|\Delta f\|_2 \right).
$$
The factor $n^{d/2}$ could be avoided if estimates were made in $\|\cdot \|_1$.
For our goals it is much easier to work with the second norm.
Fortunately, numerical experiments for the Newton potential from the following section illustrates that this factor is overestimated in practice.
Moreover for certain kernels matrix $G$ may have bounded second norm.



\subsection{Algorithm complexity in different formats} \label{sec-complexity}
Let us estimate the complexity for different formats. For simplicity, in the complexity estimates we assume that tensors $\bold{c}_g$, $\bold{q}_f$ and $\bold{\tilde w}$ 
have $n_k \sim n$ and $r_k \sim r$.  Our additional assumption in the complexity estimates is that the result of the convolution \emph{can be well approximated with the ranks} $R_k \ll r^2$.
This assumption has to be verified for each particular case, but it is standard for such kind of algorithms. 
\paragraph{Skeleton decomposition}
First, consider two-dimensional case. In two dimensions the only way to separate variables  is to approximate a matrix $A\in \mathbb{C}^{n \times m}$ by a skeleton decomposition:
\[
	A \approx U V^{T},
\]
where $U\in\mathbb{C}^{n \times r}$, $V\in\mathbb{C}^{m \times r}$ and $r$ is an approximate rank of the matrix $A$. 
Cross algorithms to compute the skeleton decomposition require $r$ columns and $r$ rows to be computed. 
Computation of a column or a row of a matrix given by its skeleton decomposition costs $\mathcal{O}(nr)$ operations.
Indeed, consider the computation of the $j$-th column:
\[	
	A(:, j) = U V(j,:)^T.
\]
The computation of the product $U V(j,:)^T$ requires $nr$ operations. As a result, the evaluation of $r$ crosses of the matrix $\bold{\Theta} = \mathcal{F} (\bold{c}_g)  \circ \mathcal{F} (\bold{q}_f)$ from step 2 needs $\mathcal{O}(nr^2)$ flop.
Additional operations performed in the cross approximation methods also have $\mathcal{O}(nr^2)$ complexity \cite{tee-cross-2000, bebe-2000}. Note, that the FFT operations from steps 1 and 3 cost $\mathcal{O}(rn \log n)$ operations.
Thus,  the algorithm complexity in the two-dimensional case is $\mathcal{O}(nr^2 + rn \log n)$.

\paragraph{Tucker format}
The Tucker format contains exponential in $d$ number of parameters $\mathcal{O}(r^d + nrd)$, but it can be efficient for problems with small $d$, especially for the case $d = 3$.
Let us calculate the complexity of the three-dimensional convolution in the Tucker format. Several implementations of cross-types  methods for the Tucker format are available,
with the first one (Cross3D) proposed in \cite{ost-tucker-2008}, see also \cite{khor-ml-2009, bebe-aca-2011} for other approaches.  For the numerical experiments in this paper we implemented a new variant of the Cross3D method -- Schur-Cross3D which has better asymptotic complexity in $r$ than the method described in \cite{ost-tucker-2008}. Details of the implementation may be found in Appendix \ref{sec-cross3d}. This method requires the computation of \emph{fibers} (which are three-dimensional generalization of columns and rows). For the interpolation,
$r$ fibers in each direction must be computed.  Let us estimate the complexity of such computation, when our approximated tensor is given as an element-wise product of two 
tensors in the Tucker format.
Let tensors $\mathcal{F} (\bold{c}_g)$ and $\mathcal{F} (\bold{q}_f)$ be in the Tucker format. 
A fiber is defined by two fixed indices, for example,  let $\bold{A}$ be a tensor, 
$$
	\bold{A}(i_1, i_2, i_3) = \sum_{\alpha_1, \alpha_2, \alpha_3} G^{(\bold{A})} (\alpha_1, \alpha_2, \alpha_3) \, U_1^{(\bold{A})}(i_1, \alpha_1) \, U_2^{(\bold{A})}(i_2, \alpha_2) \, U_3^{(\bold{A})}(i_3, \alpha_3),
$$
and $i_2, i_3$ are the fixed indices (and we need to compute the result for all $i_1 = 0, \ldots, n_1-1$).
First, we calculate 
$$
	B_{i_2 i_3}(\alpha_1) = \sum_{\alpha_2, \alpha_3} G^{ (\bold{A})} (\alpha_1, \alpha_2, \alpha_3) \,  U_2^{(\bold{A})}(i_2, \alpha_2) \, U_3^{(\bold{A})}(i_3, \alpha_3),
$$
and this step requires $\mathcal{O}(r^3)$ operations. 
Then a first mode fiber is
$$
	\bold{A}(:, i_2, i_3) = \sum_{\alpha_1} U_1^{(\bold{A})}(:, \alpha_1) B_{i_2 i_3}(\alpha_1), 
$$
and that step requires $\mathcal{O}(nr)$ operations. 
Thus, the computation of  one fiber of $\bold{\Theta} = \mathcal{F} (\bold{c}_g)  \circ \mathcal{F} (\bold{q}_f)$ costs $\mathcal{O}(nr + r^3)$ flop.
Since the Schur-Cross3D method uses $r$ fibers in each direction, the element-wise product complexity is $\mathcal{O}(nr^2 + r^4)$.

As in two-dimensional case, for a tensor $\bold{A}$ in the Tucker format, the Fourier transform $\mathcal{F} (\bold{A})$ does not change its Tucker ranks and is equivalent to three one-dimensional FFTs of the Tucker factors. 
Thus, the complexity of the steps 1 and 3 is $\mathcal{O}(nr \log n)$. The total complexity for the approximate convolution in the Tucker format is  $\mathcal{O}(nr^2 + rn \log n + r^4)$ flop.

\paragraph{TT format}
For high dimensions the Tucker format becomes unusable, and the TT-format or HT-format that have linear scaling with $d$ should be used.
A cross method for the TT-format was proposed first in \cite{ot-ttcross-2010}  and later significantly improved in \cite{so-dmrgi-2011proc} and \cite{sav-qott-2014} (and possible 
improvements are still on the way!). The asymptotic complexity of those algorithms in our case can be shown  to be equal to $\mathcal{O}(dnr^3)$ flop. 
The algorithm consists in $d$ multiplications of matrices of size $r \times r$ by matrices of size $r \times nr$. 
The FFT step can be implemented via one-dimensional FFTs of each TT-core and it costs $\mathcal{O}(r^2n \log n)$ flop.
The final algorithm complexity is $\mathcal{O}(dnr^3 + r^2 n \log n)$ and that possibly allows for very large $n$ and $d$.
\paragraph{HT and extended TT formats}
If $n$ is very large, additional complexity reduction can be achieved by using either HT or extended TT-formats.
A variant of the cross method for the HT-format can be found in \cite{lars-htcross-2013}. HT-format can be considered as a sequential application of the Tucker decomposition,
while the extended TT-format uses one preliminary Tucker decomposition and then applies TT-decomposition to the Tucker core. 
Note that the TT-format can be considered as a special case of the HT-format with a linear reduction tree. However, there is a freedom in the TT-format since it is different for different orderings of the indices. In practice TT-format is often found to be much more simple to work with, however there are examples of tensors where the HT-format gives better approximation \cite{grhack-introcmam-2011}.
If we assume that all ranks are bounded by $r$, the complexity of the convolution algorithm will be $\mathcal{O}(dnr^2 + dr^4+ rn \log n)$ flop. Note that it has better complexity 
with respect to $n$. 
The complexity estimates are summarized in the Table \ref{table:complexity}.
\begin{table}[h]
\begin{center}
\caption{Cross-conv complexity in different formats}
\label{table:complexity}
\begin{tabular}[t]{ll}
\hline
Format & Complexity \\
\hline
Skeleton decomposition & $\mathcal{O}(nr^2 + rn \log n)$ \\
Tucker 3D & $\mathcal{O}(nr^2 + r^4 + rn\log n)$\\
TT & $\mathcal{O}(dnr^3 + r^2n \log n)$\\
HT/extended TT & $\mathcal{O}(dnr^2 + dr^4 + rn\log n)$\\
\end{tabular}
\end{center}
\end{table}

\section{Numerical experiments}\label{sec-experiments}

In the numerical experiments we consider a three-dimensional case and the Tucker format. We use a new implementation of the Cross3D approximation algorithm -- Schur-Cross3D.
Schur-Cross3D and cross-conv algorithms are implemented in Python. Their implementation and the toolbox of basic tensor operations can be found at \url{https://github.com/rakhuba/tucker3d}.  
The version of numerical experiments described in this paper can be found at \url{https://bitbucket.org/rakhuba/crossconv-experiment}. Molecule data is provided as well.
For the basic linear algebra tasks the MKL library is used. Python and MKL are from the Enthought Python Distribution (EPD 7.3-1, 64-bit) \url{https://www.enthought.com}. Python version is 2.7.3. MKL version is 10.3-1.
Tests were performed on 4 Intel Core i7 2.6 GHz processor with 8GB of RAM. However, only 2 threads were used (this is default number of threads for MKL).
We would like to emphasize that implementation of the whole algorithm is in Python and time performance can be considerably improved by implementing the most time-consuming parts of it in C or Fortran languages.
\subsection{Newton potential in 3D}
As the first example we consider a computation of the Newton potential which is the convolution with $1/r$ in three dimensions
\begin{equation}\label{3d-pot}
	V(x)= \left( f*\frac{1}{\|\cdot\|}\right)(x) \equiv \int_{\mathbb{R}^3} \frac{f(y)}{\| x- y \|} \, dy,
\end{equation}
where $x = (x_1, x_2, x_3)\in \mathbb{R}^3$ and $\|x\| = \sqrt{x_1^2+x_2^2 + x_3^2}$.
Convolutions of such type typically arise in electronic structure calculations and serve as a testbed for different low-rank methods.
To discretize \eqref{3d-pot} we use a  Nystr\"om-type scheme \eqref{nystrom} on two shifted uniform $n\times n \times n$ grids.
The discretization error can be shown to be $\mathcal{O}(h^2 \left| \log h \right|)$ where $h$ is the mesh size.
\paragraph{Comparison with QTT matrix-by-vector multiplication}
First we compare our algorithm  with the algorithm of \cite{khkaz-conv-2013}, based on the matrix-by-vector multiplication in the QTT-format (further QTT algorithm).
We used the MATLAB implementation that is available as a part of the TT-Toolbox \cite{tt-toolbox}, 
and also replaced the fast DMRG approximate matrix-by-vector product \cite{Os-mvk2-2011} used in the original article  by a more efficient  AMEN-based matrix-by-vector product \cite{ds-amen-2014}.
The complexity of the QTT algorithm is logarithmic in the mode size. However, QTT ranks may be considerably larger than the Tucker ranks. 
Therefore, there is a mode size interval where the cross-conv algorithm is faster despite the fact it is asymptotically slower. 
\begin{figure}[h!]
\hfill
\begin{center} 
\resizebox{0.8\textwidth}{!}{\input{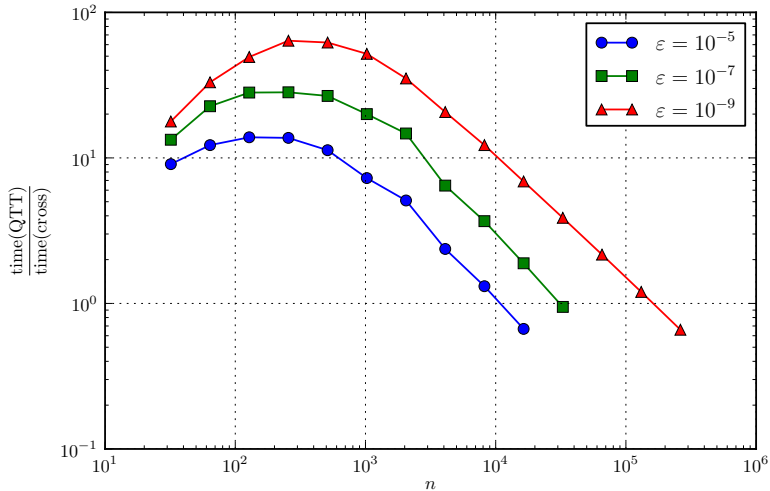}}
\end{center}
\hfill
\caption{Ratio of QTT algorithm time to cross algorithm time as a function of $n$} 
\label{pic:3d}
 \end{figure}
To illustrate this fact we consider the computation of the Newton potential of a Slater function $f(y) = e^{-\zeta|y|}$ with $\zeta =1$. Figure \ref{pic:3d} shows the ratio of the computational times as a function of the mode size. The actual timings are given in Table \ref{table:times3d}. It is clear that the more accurate the computations are (or the bigger ranks are), the faster the cross algorithm is with respect to the QTT algorithm. Moreover, it is always faster in a practically interesting range $n \sim 10^3 - 10^4$.

\begin{table}[h]
\caption{Newton potential of a Slater function}
\label{table:times3d}
\renewcommand\tabcolsep{3pt}
\begin{center}
    \begin{tabular}[H]{r|cccccccccccc}
n / $\epsilon$ &    $2^7$ & $2^8$ & $2^9$ & $2^{10}$ & $2^{11}$ & $2^{12}$ & $2^{13}$ & $2^{14}$ & $2^{15}$ & $2^{16}$ & $2^{17}$ & $2^{18}$  \\
\hline
\multicolumn{13}{c}{Cross-conv timing (sec)}\\
\hline
$10^{-5}$&   0.03 & 0.04 &0.063 & 0.12 & 0.2 & 0.5 & 1.0 & 2.2 & \\
$10^{-7}$&  0.061 & 0.091 & 0.13 & 0.24 & 0.41 & 1.1 & 2.3 & 5.2 & 11.5 \\
$10^{-9}$&   0.14 & 0.19 & 0.3 & 0.5 & 0.96 & 2.0 & 4.1 & 8.6 & 17.5 & 35.3 & 70.9 & 142.3  \\
\hline
\multicolumn{13}{c}{QTT timing (sec)}\\
\hline
$10^{-5}$&   0.42 & 0.56 & 0.71 & 0.87 & 1.0 & 1.2 & 1.3 & 1.5 & \\
$10^{-7}$&   1.7 & 2.5 & 3.6 & 4.8 & 6.0 &  7.2 & 8.6 & 9.8 & 10.9  \\
$10^{-9}$&   7.1 & 12.1 & 18.4 & 25.9 & 34.0 & 41.9 & 50.3 & 59.0 & 67.7 & 76.4 & 85.1 & 93.8  \\
\hline
\multicolumn{13}{c}{3D FFT timing (sec)}\\
\hline
 &   1.3 &  12.6 & 118.7 & 1120 &  3 hours &  &  & &
\end{tabular}
\end{center}
\end{table}

\paragraph{Newton potential of different molecules}
One of the applications of fast computation of the Newton potential are the electronic structure computations, where  the function $\rho$ is the \emph{electron density}.
We take precomputed values of $\rho$ in the Tucker format. The one-dimensional mode size is $n=5121$. The Tucker format representation was computed in \cite{st-chem-2009} and the
data was  kindly provided to us by Dr. Dmitry Savostyanov. 
Convolution times for different molecules are presented in Table \ref{tabular:molec}.
\begin{table}[h!]
\caption{Timing for the computation of the Newton potential for different molecules on the $n^3 = 5121^3$ grid}
\label{tabular:molec}
\begin{center}
\begin{tabular}[t]{lcccc}
\hline
Molecule & Accuracy & $\bold{q}_f$ ranks & $\bold{\tilde w}$ ranks & Time (s)\\
\hline
$\text{CH}_4$ &  $10^{-5}$ &  $26 \times 26 \times26$ & $22 \times22 \times22$  & 1.3 \\
		        &  $10^{-7}$ &   $39\times 39 \times39$&   $39\times 39 \times39$ & 4.1  \\
		        &  $10^{-9}$ &   $52 \times52 \times52$&   $58 \times58\times 58$&  6.4 \\
\hline
$\text{C}_2\text{H}_6$ &  $10^{-5}$ & $19 \times30\times 27$  &  $15 \times23 \times20$ &  1.2\\
		        	     &  $10^{-7}$ &  $28\times 49 \times40$ &  $24\times 42\times 39$ & 3.9 \\
		        	     &  $10^{-9}$ &  $42\times 66\times 57$ &  $39\times 66\times 60$ &  6.2\\
\hline
$\text{C}_2\text{H}_5 \text{OH}$ &  $10^{-5}$ & $43 \times42 \times43$  &  $28 \times28 \times29$ & 2.3 \\
		        			     &  $10^{-7}$ & $66 \times67\times 69$  & $50\times 50\times 51$  &  7.5\\
		        			     &  $10^{-9}$ &  $91\times 90\times 94$ &  $78\times 79\times 81$ & 19.8 \\
\hline
$\text{C}_2 \text{H}_5 \text{NO}_2$ &  $10^{-5}$ &  $24\times 60\times 60$ &  $15\times 33\times 33$ & 2.6 \\
		        			          &  $10^{-7}$ & $35\times 93 \times96$  & $26 \times61\times 62$  &  9.4\\
		        				  &  $10^{-9}$ & $45 \times126 \times133$  &  $42\times 97\times 100$ & 18.4 \\
\end{tabular}
\end{center}
\end{table}

The local filtration algorithm used in \cite{st-chem-2009} has formal complexity $\mathcal{O}(nr^2+r^5)$ for the convolution of a tensor in the canonical format 
with a tensor in the Tucker format. For the convolution of two Tucker tensors its complexity is  $\mathcal{O}(nr^2+r^6)$ compared with $\mathcal{O}(nr^2+r^4)$ complexity of the cross-conv.
We did our own implementation of the Tucker-Tucker case from \cite{st-chem-2009} and found that the Tucker ranks after local filtration are not small. For instance, Tucker ranks of the Newton potential of $\text{C}_2\text{H}_6$ are $361\times 589\times 532$ before the filtration and $82\times144\times140$ after the filtration,  while the actual ranks are $19 \times31\times 28$. Due to the strong rank dependence this leads to significantly larger computational 
time: thus, the cross-conv algorithm is more robust than the local filtration algorithm.



\begin{subsection}{Hartree-Fock equation for systems with one closed shell}
Three-dimensional convolution appears a substep in the solution of Hartree-Fock or Kohn-Sham equations in electronic structure computations. This is a classic topic and a lot of software packages are available. Here we report only preliminary results. In a series of papers Khoromskij and Khoromskaia have first used grid-based tensor methods for the solution of the Hartree-Fock equation \cite{vkhs-2el-2013,vekh-hartree-2008,mpi-chem3d-2009,venera-phd,khor-chem-2011,vekh-qtt-hartree-2011,vekh-blackbox-2013}. However, the methods they proposed are not fully ``black-box'', since they still require storage of the global basis functions for the solution, and that introduces a basis set error into the solution. We would like to store the solution of the Hartree-Fock equation as a function on a grid, i.e., as a tensor. Good news is that the solution process can be implemented solely in terms of convolutions. This is a topic of ongoing work, and in this paper we present a preliminary numerical experiment for the simplest possible case. A closed-shell Hartree-Fock equation for atoms or molecules with 2 electrons has the form
\begin{equation}\notag
	 \left(-\frac{1}{2} \,\Delta + V \right)  \psi = E \psi, \text{ where } V(x) = - \sum_{\alpha} \frac{Z_\alpha}{\|x-R_\alpha\|} + \int_{\mathbb{R}^3} \frac{|\psi(y)|^2}{\| x- y \|},
\end{equation}
$Z_\alpha$ and $R_\alpha$ are charges and coordinates of nuclei, $\psi(x)$ is the only unknown spatial orbital with $x\in\mathbb{R}^3$ and $E$ is the Hartree-Fock energy. 
Instead of the classic self consistent field (SCF) iterations we use the integral iterations (see \cite{beylkin-quantum-2004})
$$
\hat\psi = -2 (-\Delta - 2E)^{-1} V\psi \equiv -2 \, (V\psi)*\frac{e^{-\sqrt{-2E}\|\cdot\|}}{4\pi\|\cdot\|},
$$ 
(with $\hat\psi$ normalized after each iteration step) where $E$ is also recomputed at each iteration step as
$
\hat E = E + {(\hat\psi, V\hat\psi - V\psi) }/{\|\hat\psi\|^2}.
$
Note that $(-\Delta - 2E)^{-1}$ is an integral operator which is computed via the convolution with Yukawa kernel. 
At each iteration arising convolutions with Newton and Yukawa kernels are discretized via the symmetric Galerkin scheme \eqref{galerk} with piecewise-constant basis functions.
Note again that  the grid-based HF solver does not suffer from the basis set error and one can achieve necessary precision by taking larger and larger grids, i.e. reach the \emph{Hartree-Fock limit}. 
Table \ref{hf limit} illustrates this fact for the Helium atom. 
The value of the HF-limit was taken from \cite{nist-2001}. 
\begin{table}[H]
\caption{Helium atom. Dependence of the error in the Hartree-Fock energy from the grid size. $\epsilon = 10^{-6}$}
\label{hf limit}
\begin{center}
\begin{tabular}[t]{c|ccccc}
\hline
$n$ &    $1024$ & $2048$ & $4096$  & Extrapolation & HF limit (E)\\
\hline
$E_h$, (Hartree) &   -2.86113 & -2.86152 & -2.86164 & -2.861682 & -2.861679 \\
$\frac{|E_h - E|}{E}$ &  1.9e-04&   5.3e-05  & 1.25e-05 & 0.96e-06 & -\\
Time, (s) &   3.8 & 7.9& 14.3  & - & -\\
\end{tabular}
\end{center}
\end{table}
\end{subsection}
\section{Conclusion and future work}
We have presented a new efficient \emph{cross-conv} algorithm for the approximate computation of multidimensional convolution in low-rank tensor formats. The numerical experiments show that it is more efficient than the recently proposed QTT approach in a range of practically interesting mode sizes (up to $n\sim 10^4$), and the gain is higher for higher approximation accuracies or ranks. Further research will include applications of the cross-conv algorithm to a number of practically interesting models such as the Hartree-Fock and Smoluchowski equations.
\section{Acknowledgements}
We thank DrSci. Boris Khoromskij and Dr. Venera Khoromskaia for useful comments on the draft of the manuscript, and Dr. Dmitry Savostyanov for providing the data 
for the molecular densities.
We also thank anonymous referees for their comments and constructive suggestions. 
\appendix 
\section{Schur-Cross3D} \label{sec-cross3d}
Recall that in three dimensions Tucker decomposition contains only $r^3 + 3nr$ parameters and it is a natural question if it is feasible to construct this decomposition without calculating the whole three-dimensional array with $n^3$ parameters.
The answer to this question was first proposed in \cite{ost-tucker-2008}. 
This algorithm utilizes  $\mathcal{O}(nr)$ elements and has $\mathcal{O}(nr^3)$ complexity. 
It leads to $\mathcal{O}(nr^3 + r^4 + rn\log n)$ complexity for the cross-conv algorithm (see section \ref{sec-complexity} for details) which is not as efficient as the convolution based on the idea of local filtration \cite{st-chem-2009}. 
To be faster we propose a new implementation of the Cross3D method -- Schur-Cross3D with $\mathcal{O}(nr^2 + r^4)$ complexity.
Given the tensor $\bold{A}$ has exact ranks $r$, its Tucker decomposition may be represented using only $r^3 + 3nr$ elements of $\bold{A}$ as follows \cite{caiafa-cross-2010}
\begin{equation}\label{skel3d}
\bold{A} = \bold{\hat A}\times_1 U_1 \hat U_1^{-1} \times_2 U_2 \hat U_2^{-1} \times_3 U_3 \hat U_3^{-1}
\end{equation}
where $U_i$, $i=1,2,3$ consist of $r$ linearly independent fibers of the corresponding unfoldings of $\bold{A}$, $\mathcal{I}_i$, $i=1,2,3$ are numbers of $r$ linearly independent rows in matrices $U_i$,
$\bold{\hat A}= \bold{A}(\mathcal{I}_1, \mathcal{I}_2, \mathcal{I}_3)$   and 
$\hat U_i = U_i (\mathcal{I}_i, :) $.
Consider the case when tensor can be approximated with accuracy $\epsilon$ as a tensor of rank $r$.
If $\hat U_i$ is a submatrix of maximum volume in the corresponding unfolding $\bold{A}_{(i)}$ there is an approximation estimate \cite{gor-cross-2008}.
The problem to use \eqref{skel3d} with $\hat U_i$ of maximum volume is that finding the maximum volume submatrix is an NP-hard problem.
Fortunately one can use ``greedy'' strategy which is called $\verb|maxvol|$ algorithm \cite{gostz-maxvol-2010} and find quasi-maximum volume submatrices.
As a result we will get \eqref{skel3d} representation, but possibly with overestimated ranks $cr$, where $c\gtrsim1$. Numerical experiments on examples described in the previous section showed that $c\approx 1.3$.
\subsection{Theoretical estimates}
Theoretical estimates for the cross approximation are a tricky issue. 
Since cross approximation does not sample the full matrix (or tensor) it is very easy to come up 
with artificial counterexamples for any sampling technique. However, in practically interesting examples 
the convergence is very good, and that means that such matrices and tensors come from a ``good'' subclass. 
The constructive description of this subclass is an unsolved problem, however, there are several important 
theoretical results that should be mentioned. 

If the matrix (tensor) is exactly low-rank, 
then the skeleton decomposition is exact, and if during the sampling procedure we do not 
encounter zero fibers, the procedure is guaranteed to converge. In the approximate low-rank case, 
the error is multiplied by some factor. The maximum-volume principle \cite{gt-maxvol-2001} state
that if the selected rows and columns contain \emph{maximum volume} submatrix, then the error can be
estimated as
$$
   \Vert A - A_{\mathrm{skel}} \Vert_C \leq (r+1) \sigma_{r+1}.
$$
 This result was generalized to three-dimensional
and multidimensional cases in \cite{ost-tucker-2008,gor-cross-2008,sav-qott-2014}.
In practice certain greedy methods are used. For the adaptive cross approximation of 
\emph{function-generated matrices} the convergence estimate was obtained in \cite{bebe-2000,tee-cross-2000}.
An important result was obtained recently in \cite{cjd-randcrs-2013} for a class of matrices
of the form 
$$
 A = U \Phi V^{\top},
$$
where $U \in \mathbb{R}^{n \times r}$, $V \in \mathbb{R}^{m \times r}$ are orthonormal matrices and
$\Phi$ is an $r \times r$ matrix and $U$ and $V$ are $\mu$-coherent (i.e., $\max_{ij} | U_{ij} | \leq {\mu/\sqrt{n}}$), 
then it is sufficient to sample $l = \mathcal{O}(r \log n)$ columns to get an estimate with high probability. 
These results can be generalized to other SVD-based formats (Tucker, TT and HT formats) since such formats can be considered 
as sequential application of the SVD to auxiliary matrices.
\subsection{Algorithm description}


On the first step of the algorithm one may choose several randomly generated fibers or calculate $\verb|maxvol|$ fibers from an initial guess.
The goal of each next step is to add $r_0$ ``good'' fibers in the sense that they are linear independent enough to the previous ones.
Note that $r_0$ is a parameter of the algorithm and may influence convergence. 
We chose $r_0\sim 1 - 4$.

Let us consider the algorithm in more details.
Suppose that we are given $\bold{\hat A}^{(K-1)}$ of size $(K-1)r_0\times (K-1)r_0 \times (K-1)r_0$ and $U^{(K)}_i$, $i=1,2,3$ of size $n\times Kr_0$ where $K$ is the iteration number.
To find $\bold{\hat A}^{(K)}$ we should find the maximum volume submatrices in matrices $U^{(K)}_i$.
There is no guarantee that the $\verb|maxvol|$ submatrix in $U^{(K)}_i$ contains all rows of the $\verb|maxvol|$ submatrix in $U^{(K-1)}_i$. 
This leads to additional number of operations to recompute the subtensor $\bold{\hat A}^{(K-1)}$ and to add $Kr_0$ new columns into $U^{(K)}_i$ instead of $r_0$.
To avoid this one can find $r_0$ most linear independent rows to the rows of $U^{(K-1)}_i$.
We propose to do it via the Schur complement.
Thus we calculate 
$$i_1 = \verb|maxvol|(S_1),\quad i_2 = \verb|maxvol|(S_2),\quad i_3 = \verb|maxvol|(S_3)$$ 
-- indices of new rows and calculate 
$$\bold{\hat A}^{(K)} \equiv \bold{A}(\mathcal{I}^{(K)}_1, \mathcal{I}^{(K)}_2, \mathcal{I}^{(K)}_3)$$ where $$\mathcal{I}^{(K)}_1 = \mathcal{I}^{(K-1)}_1\cup i_1,\quad \mathcal{I}^{(K)}_2 = \mathcal{I}^{(K-1)}_2\cup i_2,\quad \mathcal{I}^{(K)}_3 = \mathcal{I}^{(K-1)}_3\cup i_3$$.

Next step is to find ``good'' fibers of $\bold{\hat A}^{(K)}$ in each direction to add them into matrices $U^{(K)}_i$.
To do so we calculate unfoldings $\bold{A}_{(1)}$, $\bold{A}_{(2)}$, $\bold{A}_{(3)}$ and find new ``good'' rows via the $\verb|maxvol|$ in the Schur complement of unfoldings as was mentioned before.

Finally we add $\verb|maxvol|$ fibers from the corresponding unfolding into $U^{(K)}_i$, $i=1,2,3$ and calculate $$U^{(K+1)}_i \left(\hat U^{(K+1)}_i\right)^{-1}$$ via Algorithm \ref{alg-update}. 
This algorithm allows to find new ``good'' rows without changing ``good'' rows from previous iterations.
Assuming that $r_0 \ll r \ll n$ total complexity of Algorithm \ref{alg-update} is approximately $(2r_0+1)nr$.


\begin{algorithm}[H] 
\begin{algorithmic}[1]
\caption{Factor update by Schur complement}\label{alg-update}
\REQUIRE $U\in\mathbb{C}^{n\times r}$, $u\in \mathbb{C}^{n\times r_0}$, $r_0\leqslant r$ and $\mathcal{U} = U \left[U(\mathcal{I},:) \right]^{-1}$, where $\mathcal{I}$ is a multi-index of size $r$
\ENSURE  $\mathcal{U}^{new} = U^\text{new} \left[ U^\text{new} (\mathcal{I}^{new},:) \right]^{-1}$, where $U^\text{new} = [U\ | u]$, $\mathcal{I}^{new} =  \mathcal{I}\cup i_0$ and $i_0$ is a multi-index of size $r_0$
\STATE $S = u - \mathcal{U}  u(\mathcal{I},:)$ \hfill $nr_0 + nr_0 r$
\STATE $i_0 = \verb|maxvol|(S)$ \hfill $nr_0^2$
\STATE $U_2 = S \left[S(i_0,:)\right]^{-1}$ \hfill $nr_0^2 + r_0^3$
\STATE $U_1 = \mathcal{U} - U_2  \mathcal{U}(i_0,:)$\hfill $nr + r_0 rn$
\STATE $\mathcal{U}^{new}= [U_1\ |U_2]$ 
\end{algorithmic}
\end{algorithm}

Step-by-step Schur-Cross3D is presented in Algorithm \ref{alg-cross3d}. 
Note that there are different ways to measure the accuracy of approximation. 
For instance, one may compare new good fibers with their approximation on the current iteration.

 \begin{algorithm}[H] 
\begin{algorithmic}[1]
\caption{Schur-Cross3D}\label{alg-cross3d}
\REQUIRE Function $A(i,j,k)$ which calculates certrain element of tensor $\bold{A}$, accuracy $\epsilon$ and $r_0$ -- number of fibers to be added on each iteration
\ENSURE Tucker decomposition of $\bold{A}$:  $\bold{A} \approx \bold{\hat A} \times_1 \mathcal{U}_1 \times_2 \mathcal{U}_2\times_3 \mathcal{U}_3+ \bold{E}$, $\|\bold{E}\|\leqslant \epsilon$
\STATE choose indices $\mathcal{I}_1,\mathcal{I}_2,\mathcal{I}_3$ from initial approximation or randomly
\WHILE{error $>\epsilon$}
\STATE update $\bold{\hat A} = \bold{A}(\mathcal{I}_1,\mathcal{I}_2,\mathcal{I}_3)$ 
\STATE calculate $r\times r^2$ unfoldings of $\bold{\hat A}$: $\bold{\hat A}_{(1)}, \bold{\hat A}_{(2)}, \bold{\hat A}_{(3)}$
\STATE using Schur complement find indices of new fibers $u_1,u_2,u_3$ in unfoldings $\bold{\hat A}_{(1)}, \bold{\hat A}_{(2)}, \bold{\hat A}_{(3)}$ \hfill $3r_0r^3$
\STATE error may be estimated as a norm of difference between $u_1, u_2, u_3$ and their approximation given by $\bold{\hat A},\mathcal{U}_1, \mathcal{U}_2, \mathcal{U}_3$ on the current iteration \hfill $3r_0(nr + r^3)$
\STATE add $u_1,u_2,u_3$ in $\mathcal{U}_1,\mathcal{U}_2,\mathcal{U}_3$  and find $i_0^1,i_0^2,i_0^3$ -- indices of new ``good'' rows via Algorithm \ref{alg-update} \hfill $3(2r_0 + 1) nr$
\STATE set $\mathcal{I}_1 := \mathcal{I}_1 \cup i_0^1$, $\mathcal{I}_2 := \mathcal{I}_2 \cup i_0^2$, $\mathcal{I}_3 := \mathcal{I}_3 \cup i_0^3$
\STATE set $r := r+r_0$
\ENDWHILE
\end{algorithmic}
\end{algorithm}
Thus, the overall complexity of Schur-Cross3D is approximately
$$
3(3r_0 + 1) n\sum_{k = 1}^r k + 6r_0\sum_{k = 1}^r k^3 \approx 3(3r_0 + 1) n r^2 + \frac{3r_0}{2}r^4 = \mathcal{O}(nr^2 + r^4)
$$
Note that the algorithm additionally requires $nr + r^3$ function evaluations. 
We also provide open source implementation of the proposed algorithm at \url{https://github.com/rakhuba/tucker3d} ($\verb|multifun|$ function).

\bibliographystyle{siam}
\bibliography{bibtex/tensor,bibtex/our,bibtex/dmrg,conv}
\end{document}